\documentclass{amsart}
\usepackage{amsfonts}
\usepackage{amssymb}
\usepackage[dvips]{graphics}
\usepackage{epsfig}
\usepackage{euscript}
\usepackage{color}
\newtheorem{theorem}{Theorem}[section]
          
\newtheorem{lemma}[theorem]{Lemma}

\newtheorem{definition}[theorem]{Definition}

\newcommand{\rem}[1]{{\bf Remark:}}

\newcommand{\argmin}{\mbox{argmin}}
\renewenvironment{proof}{\noindent {\bf Proof: }}{\QED\medskip}
\def\QED{{\hspace*{\fill}{\vrule height 1ex width 1ex }\quad} 
    \vskip 0pt plus20pt}
\newcommand{\be}{\begin{equation}}
\newcommand{\ee}{\end{equation}}
\newcommand{\bea}{\begin{eqnarray}}
\newcommand{\eea}{\end{eqnarray}}
\newcommand{\beann}{\begin{eqnarray*}}
\newcommand{\eeann}{\end{eqnarray*}}
\newcommand{\vt}{\vert}
\newcommand{\V}{\Vert}
\newcommand{\nn}{\nonumber}

\begin{document}

\begin{center}
{\LARGE \bf Local Stable Manifold for the Bidirectional Discrete-Time Dynamics\\[27pt]}
{\large \bf Carmeliza Navasca\\[10pt]}
{Department of Mathematics, UCLA, Los Angeles, CA 90095-1555\\[9pt]}
{navasca@math.ucla.edu\\[9pt]}
\end{center}


\noindent
{\bf Abstract.}
We show the existence of a local stable manifold for a bidirectional 
discrete-time nondiffeomorphic nonlinear Hamiltonian dynamics.  
This is the case where zero is a closed loop eigenvalue and therefore 
the Hamiltonian matrix is not invertible.  In addition, we show the 
eigenstructure and the symplectic properties of the mixed direction 
nonlinear Hamiltonian dynamics.  We extend the Local Stable Manifold 
Theorem for the nonlinear discrete-time Hamiltonian map with a 
hyperbolic fixed point.  As a consequence, we show the existence of a 
local solution to the Dynamic Programming Equations, the equations 
corresponding to the discrete-time optimal control problem.

\vspace{8pt}
\vskip .2 cm
\noindent
{\small \bf AMS numbers: 49J99, 93C55, 34K19, 93C10, 49L99}
\vfill
\hrule width2truein \smallskip {\baselineskip=10pt \noindent Research supported in part by NSF DMS-0204390 and AFOSR F49620-01-1-0202.\par }

\newpage

\section{Introduction}
One application of the local stable manifold theorem in optimal control problems is that the stable manifold of the associated Hamiltonian dynamics describes the graph of the optimal cost.  The optimal cost and the optimal control satisfy the Dynamic Programming Equations (DPE), which are obtained from the optimal control problem of minimizing a discrete-time, nonlinear cost subject to a nonlinear discrete-time dynamics through the dynamic programming technique.  In proving the existence of the solutions of the DPE, we will use the Pontryagin Maximum Principle (PMP) which gives the bidirectional nonlinear Hamiltonian dynamics and the condition for a control to be optimal satisfied by the optimal state and costate trajectories.

In the case in which the Hamiltonian matrix of the dynamics is invertible, the nonlinear Hamiltonian dynamics can be rewritten as a dynamics with both state and costate dynamics propagating in the direction where time approaches infinity.  Assuming the invertibility of the Hamiltonian matrix is to exclude zero as a closed loop eigenvalue.  The formal solutions to the DPE have been worked out using Al'brecht's method (\cite{Al61}) in \cite{Na02} and are valid for all closed loop eigenvalues lying inside the unit circle.  To our knowledge the local stable manifold theorem has only been proven for invertible maps, see (\cite{GH86}, \cite{Ha82}).  Our main result is the extension of the local stable manifold theorem to the bidirectional, discrete-time, nondiffeomorphic, nonlinear Hamiltonian dynamics.  This generalizes the proof of the existence of the local solutions to the DPE found in \cite{Na02}.

For invertible maps, Hartmann \cite{Ha82} had shown the existence of a stable manifold by the method of successive approximations on the implicit functional equation.  Another method, developed by Kelley \cite{Ke66}, is the technique of using the Contraction Mapping Theorem on a complete space.  There is another method by Irwin \cite{Ir70} based on an application of the inverse function theorem on a Banach space of sequences.  After a two-step process of diagonalizing the bidirectional discrete-time Hamiltonian dynamics, we apply the technique of Kelley \cite{Ke66} on a complete space of Lipschitz functions endowed with the supremum norm.  

The paper is organized as follows.  In the next section we introduce the bidirectional discrete-time Hamiltonian dynamics from the Pontryagin Maximum Principle that is associated with the optimal control problem.  In Section 3 we give some discussion of Gronwall's inequalities in the discrete-time case, which will then be used in the proof.  In Section 4 we state and prove the local stable manifold theorem.  In Section 5 we discuss the eigenstructure and the symplectic properties of the dynamics.  Finally, in Section 6 we show how the Theorem along with these properties give the local solvability of the Dynamic Programming Equations.

\section{Nonlinear Dynamics}
We formulate a discrete in time infinite horizon optimal control problem of minimizing the cost functional, 
\[ \min_u \sum_{k=0}^{\infty} l(x_k,u_k)\]
subject to the dynamics
\beann
x^+&=&f(x,u)\\
x(0)&=&x_0
\eeann
where the state vector $x\in \mathbb{R}^n$, the control $u\in \mathbb{R}^m$, and 
\bea
f(x,u)&=&Ax + Bu + f^{[2]}(x,u) + f^{[3]}(x,u) + \ldots \label{f}\\
l(x,u)&=&\frac{1}{2}x'Qx + x'Su + \frac{1}{2}u'Ru + l^{[3]}(x,u) + \ldots \label{l}
\eea
where $f^{[m]}(x,u)$ and $l^{[m]}(x,u)$ as homogeneous polynomials in $x$ and $u$ of degree $m$.
We let $x^+=x_{k+1}$ and $x=x_k$.

Associated with optimal control problem formulation there is a nonlinear Hamiltonian,
\bea \label {theham}
H(x,u,\lambda^+)=(\lambda^{+})'f(x,u) + l(x,u)
\eea
where $\lambda^+=\lambda_{k+1}$ and the functions $f$ and $g$ are given by equations (\ref{f}) and (\ref{l}).
The Pontryagin Maximum Principle (PMP) states the following:
\begin{theorem}
If $x_k$ and $u_k$ are optimal for $k \in 0, 1, 2, \ldots$, then there
exists $\lambda_k \neq 0$ for $k \in 0, 1, 2, \ldots$ such that
\bea
x^+&=&\frac{\partial H}{\partial \lambda^+}(x,u,\lambda^+) \label{hode1} \\
\lambda&=& \frac{\partial H}{\partial x}(x,u,\lambda^+) \label{hode2}
\eea
and
\bea
u^*&=&\argmin_{u}H(x,u,\lambda^+). \label{hode3}
\eea
\end{theorem}
Thus, the minimizer of the nonlinear Hamiltonian evaluated at the optimal
$x$ and $\lambda^+$ is the optimal control $u^*$ amongst all admissible
controls $u$.  Note that $u^*(x,\lambda^+) \in \mathcal{C}^{r-1}$ since $f \in \mathcal{C}^{r-1}$ and $l \in \mathcal{C}^r$ in (\ref{theham}).  We assume that $H$ is convex in u to guarantee a unique optimal control.  With the Hamiltonian (\ref{theham}), the equations (\ref{hode1}) and (\ref{hode2}) are the following
\bea \label{nonham}
\left[ \begin{array}{c}
x^+ \\
\lambda \\
\end{array}\right]&=&
\left[\begin{array}{cc}
A - BR^{-1}S' & -BR^{-1}B'\\
Q - SR^{-1}S' & A' - SR^{-1}B'\\
\end{array}\right]
\left[\begin{array}{c}
x \\
\lambda^+ \\
\end{array}\right]+
\left[\begin{array}{c}
F(x,\lambda^+)\\
G(x,\lambda^+)\\
\end{array}\right]
\eea
where $x,\lambda \in \mathbb{R}^n$ and $F(x,\lambda^+)$ and $G(x,\lambda^+)$ contain the nonlinear terms.
Observe the opposing directions of the propagation of the state and costate dynamics in (\ref{nonham}). 

If we linearize the $2n$ dimensional difference equations (\ref{nonham}) around zero, we obtain exactly the linear Hamiltonian system,
\bea \label{hamsys}
\left[ \begin{array}{c}
x^+ \\
\lambda \\
\end{array}\right]=\mathbb{H}
\left[\begin{array}{c}
x \\
\lambda^+ \\
\end{array}\right],
\eea
where
\beann
\mathbb{H}=
\left[ \begin{array}{cc}
A - BR^{-1}S'  & -BR^{-1}B' \\
Q -SR^{-1}S'   &  A' - SR^{-1}B'\\
\end{array}\right]
\eeann
is the associated Hamiltonian matrix and the corresponding Hamiltonian is
\[H(x,\lambda^+,u)=\lambda^{+'}(Ax + Bu) + \frac{1}{2} x^{'} Q x + x^{'}S u
+ \frac{1}{2}u^{'}R u.\]

In fact, the stable subspace of $\mathbb{H}$ is described by
\[\lambda=Px\]
where the nonegative definite P satisfies the discrete-time algebraic Riccati equation (DTARE),
\bea \label{dtare}
P-A^{'}PA + (A^{'}PB + S)(B^{'}PB+R)^{-1})(A^{'}PB+S)^{'}-Q=0.
\eea
Moreover, the stable linear subspace is 
\beann
E_s=
\left[\begin{array}{c}
I \\
P \\
\end{array}\right]
\eeann  
where $E_s$ is spanned by the $n$ stable eigenvalue of the Hamiltonian matrix lying inside the unit circle.  A detailed proof is found in \cite{Na02}. Since the linear part of the nonlinear bidirectional dynamics is the linear Hamiltonian system, we have that the linear part of the local stable manifold is $\lambda=Px$.  

\section{Discrete-Time Version of Gronwall's Inequalities}

We first discuss various results based on the discrete-time version of Gronwall's inequalities as these following lemmas will be useful in the proof of the local existence of a stable manifold.

\begin{lemma} \label{gronwall1}
Suppose the sequence of scalars $\{u_j \}_{j=0}^{\infty}$ satisfies the difference inequality
\bea \label {groneqn}
u_{k+1} \leq  \delta u_k + L
\eea
where $\delta,~L \geq 0$, then
\[u_k \leq \delta^k u_0 + L\sum_{j=0}^{k-1}\delta^{k-1-j}.\]
\end{lemma}
The proof of the lemma above clearly follows from summing equation (\ref{groneqn}) from $1$ to $k$. 

\begin{lemma} \label{gronwall2}
Suppose $\{\xi_j\}_{j=0}^{\infty}$ is a sequence that satisfies
\[\vert \xi_k \vert \leq C_1 \sum_{j=0}^{k-1} \vert \xi_j \vert + C_2\]
with constants $C_1, C_2 \geq 0$.  Then 
\[\vert \xi_k \vert \leq C_2 \sum_{j=1}^{k}(1+C_1)^{j}.\]
\end{lemma}

\begin{proof}
Let $s_k=\sum_{j=0}^{k-1}\vt \xi_j \vt$.  Then, the sequence $\{s_j\}_{j=0}^{\infty}$ satisfies
\[s_{k+1} \leq (1 + C_1) s_k + C_2\]
where $C_1,C_2 \geq 0$.  By Lemma \ref{gronwall2},
\[\vert s_k \vert \leq (1+C_1)^k \vert s_0 \vert + C_2 \sum_{j=0}^{k-1}(1+C_1)^{k-1-j}.\]
It follows that
\beann
\vert \xi_k \vert &\leq& (1 + C_1) \vert s_k \vert + C_2\\
                  &\leq& (1 + C_1) \Big[(1+C_1)^k \vert s_0 \vert + C_2 \sum_{j=0}^{k-1}(1+C_1)^{k-1-j}\Big] \\
                  &\leq& C_2 \sum_{j=0}^{k-1}(1+C_1)^{k-j}\\
                  &\leq& C_2 \sum_{j=1}^{k}(1+C_1)^j\\
\eeann
since $\vert s_0\vert=0$. 
\end{proof}

\section{Local Stable Manifold Theorem for the Bidirectional Discrete-Time Dynamics}

In this section we prove the existence of a local stable manifold
\bea \label{phimanifold}
\lambda=\phi(x)
\eea
for the Hamiltonian dynamics,
\bea \label{nondiffham}
\left[ \begin{array}{c}
x^+ \\
\lambda \\
\end{array}\right]&=&
\left[\begin{array}{cc}
A & -BR^{-1}B'\\
Q & A'\\
\end{array}\right]
\left[\begin{array}{c}
x \\
\lambda^+ \\
\end{array}\right]+
\left[\begin{array}{c}
F(x,\lambda^+)\\
G(x,\lambda^+)\\
\end{array}\right]
\eea
where $x,\lambda \in \mathbb{R}^n$ and zero is an eigenvalue of A.  The nonlinear terms, $F$ and $G$, are $\mathcal{C}^r$ functions for $r \geq 1$ such that
\bea \label{condsnonl}
F(0,0)=0, && G(0,0)=0 \\
\frac{\partial F}{\partial (x,\lambda)}(0,0)=0, && \frac{\partial G}{\partial (x,\lambda)}(0,0)=0. \nonumber
\eea

Since we are locally proving the existence of the stable manifold, we just need the local behavior of the dynamics; so first we talk about cut-off functions.  The proof also requires the discussion on the stability of the nonlinear state dynamics.  So the following subsection deals the local asymptotic stability of the state dynamics. Then, we describe the diagonalization of the bidirectional Hamilton system.  Finally, we show the existence of $\lambda=\phi(x)$.

\subsection{Cut-off Functions}

First, we introduce a $\mathcal{C}^{\infty}$ cut-off function $\rho(y): \mathbb{R}^n \longrightarrow [0,1]$ such that
\beann
\rho(y) =
   \begin{cases}
    1, &\text{if $0 \leq \vert y \vert \leq 1$}\\
    0, &\text{if $\vert y \vert > 2$}
   \end{cases}
\eeann
and $0 \leq \rho(y) \leq 1$ otherwise.  Then we define the functions
\bea \label{locfuncs}
F(x,\lambda^+;\epsilon)&:=& F(x\rho(\frac{x}{\epsilon}),\lambda^+\rho(\frac{\lambda^+}{\epsilon}))\\
G(x,\lambda^+;\epsilon)&:=& G(x\rho(\frac{x}{\epsilon}),\lambda^+\rho(\frac{\lambda^+}{\epsilon}))\nn
\eea
for $x,\lambda^+ \in \mathbb{R}^n$.  Since the $F(x,\lambda^+)$ and $G(x,\lambda^+)$ agree with $F(x,\lambda^+;\epsilon)$ and $G(x,\lambda^+;\epsilon)$, respectively, for $\vt x \vt, \vt \lambda^+ \vt \leq \epsilon$, it suffices to prove the existence of a stable manifold for some $\epsilon > 0$. 

Now, we show that there exists $N_1>0$ and $N_2>0$ such that
\bea
\vt F(x,\lambda;\epsilon)-F(\tilde{x},\tilde{\lambda};\epsilon) \vt &\leq& N_1 \epsilon \Big[ \vert x- \tilde{x} \vert + \vt  \lambda -\tilde{\lambda} \vt \Big]\label{fbnd}\\
\vt G(x,\lambda;\epsilon)-G(\tilde{x},\tilde{\lambda};\epsilon) \vt &\leq& N_1 \epsilon \Big[ \vert x- \tilde{x} \vert + \vt  \lambda -\tilde{\lambda} \vt \Big]\label{gbnd}
\eea
and 
\bea
\Big \vt \frac{\partial F}{\partial (x,\lambda)}(x,\lambda;\epsilon)- \frac{\partial F}{\partial (x,\lambda)} (\tilde{x},\tilde{\lambda};\epsilon) \Big \vt &\leq&  N_2 \Big[ \vert x- \tilde{x} \vert + \vt  \lambda -\tilde{\lambda} \vt \Big] \label{fbnd}\\
\Big \vt \frac{\partial G}{\partial (x,\lambda)}(x,\lambda;\epsilon)-\frac{\partial G}{\partial (x,\lambda)}(\tilde{x},\tilde{\lambda};\epsilon) \Big \vt &\leq& N_2 \Big[ \vert x- \tilde{x} \vert + \vt  \lambda -\tilde{\lambda} \vt \Big]. \label{gbnd}
\eea
Since $\rho(y)$ and its partial derivatives are continuous functions with compact support there exists $M>0$ such that
\beann
\Big \vert \frac{\partial \rho}{\partial y}(y) \Big \vert &\leq& M\\
\Big \vert \frac{\partial^2 \rho}{\partial y^2}(y) \Big \vert &\leq& M
\eeann
for all $\lambda \in \mathbb{R}^n$.  We also choose $M>0$ large enough that
\bea 
\Big \vt \frac{\partial F}{\partial x}(x,\lambda)         \Big \vt &\leq& M \vt x \vt \label{bndonf1a}\\
\Big \vt \frac{\partial F}{\partial \lambda}(x,\lambda)   \Big \vt &\leq& M \vt \lambda \vt \label{bndonf1b}\\
\Big \vt \frac{\partial^2 F}{\partial x^i \partial \lambda^j}(x,\lambda) \Big \vt &\leq& M,~~~i,j=1,2 \label{bndonf2}
\eea
because of the condition (\ref{condsnonl}) for $\vt x \vt, \vt \lambda \vt <1$.
By the Mean Value Theorem,
\bea \label{mvt}
\vert F(x,\lambda;\epsilon) - F(\tilde{x},\tilde{\lambda};\epsilon)  \vert &\leq & \vt F(x,\lambda;\epsilon) - F(\tilde{x},\lambda;\epsilon) + F(\tilde{x},\lambda;\epsilon) - F(\tilde{x},\tilde{\lambda};\epsilon) \vt \nn \\
                                                                           &\leq & \vt F(x,\lambda;\epsilon) - F(\tilde{x},\lambda;\epsilon) \vt + \vt F(\tilde{x},\lambda;\epsilon) - F(\tilde{x},\tilde{\lambda};\epsilon)\vt \nn \\
                                                                           &\leq & \Big \vt \frac{\partial F}{\partial x}(\xi_1,\lambda;\epsilon)  \Big \vt \vt x-\tilde{x} \vt +\Big \vt \frac{\partial F}{\partial \lambda}(x,\xi_2;\epsilon)  \Big \vt \vt \lambda-\tilde{\lambda} \vt \nn 
\eea
where $\xi_1$ is between $x$ and $\tilde{x}$ and $\xi_2$ is between $\lambda$ and $\tilde{\lambda}$.  Similarly,
\beann
\Big \vert \frac{\partial F}{\partial x}(x,\lambda;\epsilon) - \frac{\partial F}{\partial x}(\tilde{x},\tilde{\lambda};\epsilon)  \Big \vert &\leq & \Big \vt  \frac{\partial^2 F}{\partial x^2} (\xi_1,\lambda)\Big \vt  \vt x-\tilde{x} \vt  + \Big \vt \frac{\partial^2 F}{\partial x \partial \lambda} (x,\xi_2) \Big \vt  \vt \lambda -\tilde{\lambda}  \vt  \\
\Big \vert \frac{\partial F}{\partial \lambda}(x,\lambda;\epsilon) - \frac{\partial F}{\partial \lambda}(\tilde{x},\tilde{\lambda};\epsilon)  \Big \vert &\leq & \Big \vt  \frac{\partial^2 F}{\partial \lambda \partial x} (\xi_1,\lambda)\Big \vt  \vt x-\tilde{x} \vt  + \Big \vt \frac{\partial^2 F}{\partial \lambda^2} (x,\xi_2) \Big \vt  \vt \lambda -\tilde{\lambda}  \vt.
\eeann

Next we estimate for $0 \leq \epsilon < \frac{1}{2}$
\beann
\Big \vt \frac{\partial F}{\partial x}(\xi_1,\xi_2;\epsilon) \Big\vt &=& \Big\vt \frac{\partial F}{\partial x}(\rho(\frac{\xi_1}{\epsilon})\xi_1,\frac{\xi_2}{\epsilon})\xi_2)\Big\vt \; \Big\vt \frac{\partial \rho}{\partial y}(\frac{\xi_1}{\epsilon})\frac{\xi_1}{\epsilon} + \rho(\frac{\xi_1}{\epsilon}) \Big\vt \nn \\
&\leq& M \vt \rho(\frac{\xi_1}{\epsilon})\vt \vt \xi_1 \vt \Big( \Big \vt \frac{\partial \rho}{\partial y}(\frac{\xi_1}{\epsilon}) \frac{\xi_1}{\epsilon} \Big \vt + \vt \rho(\frac{\xi_1}{\epsilon}) \vt \Big) \nn \\ 
&\leq& M[M+1]\epsilon
\eeann
and 
\beann
\Big \vt \frac{\partial^2 F}{\partial x^2}(\xi_1,\xi_2;\epsilon) \Big\vt &=& \Big \vt \frac{\partial^2 F}{\partial x^2}(\rho(\frac{\xi_1}{\epsilon})\xi_1,\frac{\xi_2}{\epsilon})\xi_2) \Big \vt \; \Big\vt \frac{\partial \rho}{\partial y}(\frac{\xi_1}{\epsilon})\frac{\xi_1}{\epsilon} + \rho(\frac{\xi_1}{\epsilon}) \Big\vt \\
&&+ \Big\vt \frac{\partial F}{\partial x}(\rho(\frac{\xi_1}{\epsilon})\xi_1,\frac{\xi_2}{\epsilon})\xi_2)\Big\vt \; \Big\vt \frac{\partial \rho}{\partial \xi_1}(\frac{\xi_1}{\epsilon})\frac{2}{\epsilon} + \frac{\partial^2 \rho}{\partial y^2}(\frac{\xi_1}{\epsilon}) \frac{\xi_1}{\epsilon^2} \Big\vt \\
& \leq & M(1+2M)+8M^2
\eeann
for $\vt \xi_1 \vt, \vt \xi_2 \vt \leq \epsilon$.

Let
\beann
N_1 &=& M^2[M+1]\\
N_2 &=& M(1+2M)+8M^2.
\eeann
It follows that
\bea
\Big \vert \frac{\partial F}{\partial (x,\lambda)}(\xi_1,\xi_2;\epsilon) \Big \vert &\leq& N_1\epsilon \label{origmvt1} \\
\Big \vt \frac{\partial^2 F}{\partial x^i \partial \lambda^j}(\xi_1,\xi_2;\epsilon) \Big\vt &\leq& N_2,~~~i,j=1,2 \label{origmvt2}
\eea
for $\vt \xi_1 \vt, \vt \xi_2 \vt < \epsilon$.  The inequalities above also hold for $G$ as well.

Thus, 
\bea
\vt F(x,\lambda;\epsilon)-F(\tilde{x},\tilde{\lambda};\epsilon) \vt &\leq& N_1 \epsilon \Big[ \vert x- \tilde{x} \vert + \vt  \lambda -\tilde{\lambda} \vt \Big]\label{fbnd1}\\
\vt G(x,\lambda;\epsilon)-G(\tilde{x},\tilde{\lambda};\epsilon) \vt &\leq& N_1 \epsilon \Big[ \vert x- \tilde{x} \vert + \vt  \lambda -\tilde{\lambda} \vt \Big]\label{gbnd1}
\eea
and 
\bea
\Big \vt \frac{\partial F}{\partial (x,\lambda)}(x,\lambda;\epsilon)- \frac{\partial F}{\partial (x,\lambda)} (\tilde{x},\tilde{\lambda};\epsilon) \Big \vt &\leq&  N_2 \Big[ \vert x- \tilde{x} \vert + \vt  \lambda -\tilde{\lambda} \vt \Big] \label{fbnd2}\\
\Big \vt \frac{\partial G}{\partial (x,\lambda)}(x,\lambda;\epsilon)-\frac{\partial G}{\partial (x,\lambda)}(\tilde{x},\tilde{\lambda};\epsilon) \Big \vt &\leq& N_2 \Big[ \vert x- \tilde{x} \vert + \vt  \lambda -\tilde{\lambda} \vt \Big].
 \label{gbnd2}
\eea

Henceforth we suppress $\epsilon$ and write $F(x,\lambda), G(x,\lambda)$ for $F(x,\lambda;\epsilon),~G(x,\lambda;\epsilon)$.

\subsection{Stability of the Nonlinear Dynamics}

From Section 2, we mentioned that the linear term of the stable manifold for the nonlinear bidirectional Hamiltonian dynamics is $Px$.  Then, the local stable manifold is of the form
\bea \label{phiform}
\lambda = \phi(x)=Px + \psi(x)
\eea
where $\psi(x)$ contains all the nonlinear terms.

Suppose we substitute (\ref{phiform}) into the state dynamics in (\ref{nondiffham}), then the nonlinear state dynamics becomes
\beann
(I+BR^{-1}B'P)x^+=Ax-BR^{-1}B'\psi(x^+)+ F(x,Px^+ + \psi(x^+)).
\eeann
By the Matrix Inversion Lemma (\cite{LS95}), we have that 
\[(I+BR^{-1}B'P)^{-1}=(I-B(B'PB+R)^{-1}B'P).\]
Then, it follows that 
\bea \label{nlinearxmap}
x^+&=&(A+BK)x + f_{\psi}(x,x^+)\\
x(0)&=&x_0 \nn
\eea
where $K=-(B'PB + R)^{-1}B'P$ and $f_{\psi}(x,x^+)=(I+BR^{-1}B'P)^{-1}(F(x,Px^+ + \psi(x^+)) - BR^{-1}B'\psi(x^+))$.  

The implicit equation above can be solved.  Let $\mathcal{F}:\mathcal{N}_{\epsilon}(0) \subset \mathbb{R}^{2n} \rightarrow \mathbb{R}^n$ such that 
\[\mathcal{F}(x,x^+)=x^+ - (A+BK)x-f_{\psi}(x,x^+)=0\]
for $x,x^+ \in \mathcal{N}_{\epsilon}(0)$ where $\mathcal{N}_{\epsilon}(0)$ is an open neighborhood of radius $\epsilon$ around $0$.  Then, for $0 \in \mathcal{N}_{\epsilon}(0)$ the Jacobian
\beann
\frac{\partial \mathcal{F}}{\partial x^+}(0) &=& I-\frac{\partial f_{\psi}}{\partial x^+}(0) \\
                                             &=& I-(I+BR^{-1}B'P)^{-1}\Big(\frac{\partial F}{\partial \lambda^+}(0,\phi(0))\frac{\partial \psi}{\partial x^+}(0) - BR^{-1}B'\frac{\partial \psi}{\partial x^+}(0)\Big)\\
                                             &=& I,
\eeann
because of the condition (\ref{condsnonl}) and $\psi(x^+)$ only contains nonlinear terms.  Then, by the Implicit Function Theorem there exists $\mathbb{F}(x)$ such that
\bea \label{nlinearxmap2}
x^+&=&\mathbb{F}(x)
\eea
is equivalent to the earlier state dynamics (\ref{nlinearxmap}).  Moreover,  the linear term of $\mathbb{F}(x)$ is $(A+BK)x$, i.e.;
\bea \label{consxeq2}
\mathbb{F}(x)= (A+BK)x + F_{\psi}(x)
\eea
because
\beann
\frac{\partial \mathbb{F}}{\partial x}(0)&=&- \Big(\frac{\partial \mathcal{F}}{\partial x^+}(0)\Big)^{-1}\frac{\partial \mathcal{F}}{\partial x}(0)\\
&=&-I[-(A+BK)]\\ 
&=& A+BK. 
\eeann
It follows that $F_{\psi}(x)$ contains only the nonlinear terms and thus,
\bea \label{consxeq3}
F_{\psi}(0)=0
\eea
and
\[\frac{\partial F_{\psi}}{\partial x_i}(0)=0\;\;i=1,\ldots,n.\]

The linear part of (\ref{nlinearxmap}) is
\bea \label{linearxmap}
x^+=(A+BK)x.
\eea
Since the eigenvalues of $(A+BK)$ lie strictly inside the unit circle, the term $(A+BK)^k x_0 \longrightarrow 0$ as $k \longrightarrow \infty$. Thus, the system (\ref{linearxmap}) is asymptotically stable.  Also, it implies that there exists a unique positive definite $P$ that satisfies the Lyapunov equation
\[(A+BK)'P(A+BK)-P=-I.\]

Now we show the stability of the nonlinear dynamics
\beann
x_k &=& (A+BK)x + F_{\psi}(x)\\
x(0)&=& 0. 
\eeann

We must prove that 
\beann
\lim_{x \rightarrow 0} \frac{\vt F_{\psi}(x) \vt}{\vt x \vt}=0;
\eeann
i.e., given any $\varepsilon > 0$ and any $\psi(x)$ satisfying the conditions
\bea 
\psi(0)&=&0 \label{cond1}\\
\vt \psi(x)-\psi(\bar{x})  \vt &\leq& l(\epsilon)\vt x - \bar{x}  \vt \label{cond2}
\eea
where $l(\epsilon) \longrightarrow 0$ as $\epsilon \rightarrow 0$, there exists $\delta >0$ such that
\[ \frac{\vt F_{\psi}(x) \vt}{\vt x \vt} < \varepsilon ~~~~~\text{whenever}~~ \vt x \vt < \delta.\]
We define $\psi(x)$ to be the nonlinear term of the stable manifold in (\ref{phiform}).  The conditions (\ref{cond1}) and (\ref{cond2}) will be necessary for the proof of the local stable manifold theorem.

Recall that
\beann
x^+=\mathbb{F}(x)= (A+BK)x + F_{\psi}(x). 
\eeann
Then,
\beann
0=\mathcal{F}(x,x^+)&=& x^+ - (A+BK)x-f_{\psi}(x,x^+)\\
                    &=&(A+BK)x + F_{\psi}(x)- (A+BK)x-f_{\psi}(x,(A+BK)x + F_{\psi}(x))\\
                    &=&F_{\psi}(x)-f_{\psi}(x,(A+BK)x + F_{\psi}(x)).
\eeann
It follows that
\bea \label{Fpsi}
F_{\psi}(x)&=&(I+BR^{-1}B'P)^{-1}\Big[F(x,P((A+BK)x + F_{\psi}(x)) + \psi((A+BK)x + F_{\psi}(x))) \nn \\
           &&~~- BR^{-1}B'\psi((A+BK)x + F_{\psi}(x))\Big].
\eea
Let $\mathbb{B}_1= \V (I+BR^{-1}B'P)^{-1} \V$, $\mathbb{B}_2= \V BR^{-1}B' \V$, $\mathbb{P}=\V P \V$ and $\alpha = \max_{i} \vt \lambda_i \vt$ where $\lambda_i \in \sigma(A+BK)$ and $\vt \lambda_i \vt < 1$.   We have the following from (\ref{Fpsi}),
\beann
\vt F_{\psi}(x) \vt &\leq& \mathbb{B}_1 N_1 \epsilon \left[\vt x \vt + \vt P(A+BK)x + PF_{\psi}(x) + \psi((A+BK)x + F_{\psi}(x)))\vt \right] \\
&&~~+ \mathbb{B}_1\mathbb{B}_2 \vt \psi((A+BK)x + F_{\psi}(x))) \vt.
\eeann
because of (\ref{fbnd1}).
Using the Lipschitz condition (\ref{cond2}) for $\psi(x)$,
\beann
\vt F_{\psi}(x) \vt &\leq& \left[ \mathbb{B}_1 N_1 \epsilon + \alpha(\mathbb{B}_1N_1\epsilon \V P \V + \mathbb{B}_1N_1 \epsilon l(\epsilon) + \mathbb{B}_1\mathbb{B}_2l(\epsilon)) \right] \vt x \vt \\
&&+ \left[ \mathbb{B}_1N_1\epsilon\V P \V + \mathbb{B}_1N_1 \epsilon l(\epsilon) + \mathbb{B}_1\mathbb{B}_2 l(\epsilon)\right] \vt F_{\psi}(x) \vt.
\eeann
Solving for $\vt F_{\psi}(x) \vt$,
\beann
\vt F_{\psi}(x) \vt &\leq& \frac{\mathbb{B}_1 N_1 \epsilon + \alpha(\mathbb{B}_1N_1\epsilon \V P \V + \mathbb{B}_1N_1 \epsilon l(\epsilon) + \mathbb{B}_2l(\epsilon))}{1-\left( \mathbb{B}_1N_1\epsilon + \mathbb{B}_1N_1 \epsilon l(\epsilon) + \mathbb{B}_2 l(\epsilon)\right)}.
\eeann
Let $\delta= \frac{1-\left( \mathbb{B}_1N_1\epsilon + \mathbb{B}_1N_1 \epsilon l(\epsilon) + \mathbb{B}_2 l(\epsilon)\right)}{\mathbb{B}_1 N_1 \epsilon + \alpha(\mathbb{B}_1N_1\epsilon \V P \V + \mathbb{B}_1N_1 \epsilon l(\epsilon) + \mathbb{B}_2l(\epsilon))}~\varepsilon$.  For some $\epsilon >0$ and $\varepsilon >0$, we have that $\delta > 0$.
Then,
\beann
\vt F_{\phi}(x) \vt &\leq& \frac{\mathbb{B}_1 N_1 \epsilon + \alpha(\mathbb{B}_1N_1\epsilon \V P \V + \mathbb{B}_1N_1 \epsilon l(\epsilon) + \mathbb{B}_2l(\epsilon))}{1-\left( \mathbb{B}_1N_1\epsilon + \mathbb{B}_1N_1 \epsilon l(\epsilon) + \mathbb{B}_2 l(\epsilon)\right)}~ \vt x \vt\\
                          &\leq& \frac{\mathbb{B}_1 N_1 \epsilon + \alpha(\mathbb{B}_1N_1\epsilon \V P \V + \mathbb{B}_1N_1 \epsilon l(\epsilon) + \mathbb{B}_2l(\epsilon))}{1-\left( \mathbb{B}_1N_1\epsilon + \mathbb{B}_1N_1 \epsilon l(\epsilon) + \mathbb{B}_2 l(\epsilon)\right)}~ \delta \\
                          &\leq& \varepsilon.
\eeann
Thus, 
\[F_{\phi}(x)= o(\vt x \vt).\]
Now, we use the Lyapunov argument.
Let $v(x)=x'Px$.  Then,
\beann
\Delta v(x)&=& v(x^+)-v(x)\\
           &=& {x^+}' P x^+ -x' P x\\
           &=& [(A+BK)x-F_{\psi}(x)]'P [(A+BK)x-F_{\psi}(x)]-x'Px\\
           &=& x'((A+BK)'P(A+BK)-P)x + 2x'(A+BK)'PF_{\psi}(x)\\
           &=& -\vt x \vt^2 + 2x'(A+BK)'PF_{\psi}(x).
\eeann
since
\beann
\vt F_{\psi}(x) \vt \leq \frac{1}{3p}\vt x \vt
\eeann
and
\[\vt 2x'(A+BK)'P F_{\psi}(x)\vt \leq \frac{2}{3}\vt x \vt^2  \]
for some $p>0$.
Thus,
\beann
\Delta v(x)= -\frac{\vt x \vt^2}{3}<0.
\eeann
Therefore, the nonlinear dynamics is locally asymptotically stable uniform for all $\psi \in \mathbb{X}$.

\subsection{Diagonalization of the Hamiltonian Matrix}

Recall the bidirectional nonlinear dynamics in (\ref{nondiffham}) and the condition (\ref{condsnonl}).

By substituting 
\bea \label{lasubstitute}
\lambda=Px + \psi(x)
\eea
into the state dynamics above (\ref{nondiffham}), we get
\bea \label{newstate}
x^+&=& (A+BK)x + f_{\psi}(x,x^+)
\eea
where
$f_{\psi}(x,x^+)=(I+BR^{-1}B'P)^{-1}(F(x,Px^+ + \psi(x^+)) - BR^{-1}B'\psi(x^+))$.

As we substitute (\ref{lasubstitute}) and (\ref{newstate}) into the costate dynamics in (\ref{nondiffham}), we also add $0=(BK)'\lambda^+ -(BK)'\lambda^+$.  Then, the costate dynamics becomes
\beann \label{nlinearlammap} 
\lambda=(A+BK)'\lambda^+ + \bar{Q}x + g_{\psi}(x,x^+).
\eeann
where $\bar{Q}=Q-K'B'P(A+BK)$ and $g_{\psi}(x,x^+)=  G(x,Px^+ + \psi(x^+)) + K'B'(-\psi(x^+) -Pf_{\psi}(x,x^+)$.  

Thus, the substitution of 
\[\lambda=Px + \psi(x)  \]
into the dynamics (\ref{nondiffham}) results in a new nonlinear dynamics
\bea \label{nondiffham2}
\left[ \begin{array}{c}
x^+ \\
\lambda \\
\end{array}\right]&=&
\left[\begin{array}{cc}
A+BK & 0\\
\bar{Q} & (A+BK)'\\
\end{array}\right]
\left[\begin{array}{c}
x \\
\lambda^+ \\
\end{array}\right]+
\left[\begin{array}{c}
f_{\psi}(x,x^+)\\
g_{\psi}(x,x^+)\\
\end{array}\right]
\eea
where
\beann
f_{\psi}(x,x^+)&=&(I+BR^{-1}B'P)^{-1}(F(x,\psi(x^+)) - BR^{-1}B'\psi(x^+))
\eeann
and
\beann
g_{\psi}(x,x^+)&=&G(x,\psi(x^+)) + K'B'(-\psi(x^+) -Pf_{\psi}(x,x^+).~~~~~~~~
\eeann
The nonlinear terms $f_{\psi}$ and $g_{\psi}$ are $\mathcal{C}^k$ functions for $k \geq 1$ such that
\bea 
f_{\psi}(0,0)=0, && g_{\psi}(0,0)=0 \label{fgcond1}\\
\frac{\partial f_{\psi}}{\partial(x,x^+)}(0,0)=0, && \frac{\partial g_{\psi}}{\partial (x,x^+)}(0,0)=0.\label{fgcond2}
\eea
because of (\ref{condsnonl}), $\psi(x)$ only contains nonlinear terms and 
\[\frac{\partial \psi}{\partial x}(0)=0.\]

Now we introduce the $z$ coordinate by the transformation
\bea \label{sectrans}
\lambda=z + Sx
\eea
for some matrix $S$ to block diagonalize the block lower triangular Hamiltonian matrix in (\ref{nondiffham2}).  By substitution, the system (\ref{nondiffham2}) becomes

\beann
x^+&=& (A+BK)x + f_{\psi}(x,x^+)\\
  z&=& (A+BK)'z^+ + (A+BK)'S(A+BK)x - Sx + \bar{Q}x + h_{\psi}(x,x^+)
\eeann
where 
\[h_{\psi}(x,x^+)=(A+BK)'Sf_{\psi}(x,x^+) + g_{\psi}(x,x^+).\]
From the $z$ dynamics above observe that the terms
\[(A+BK)'S(A+BK)x - Sx + \bar{Q}x=0\]
where $\bar{Q}=Q-K'B'P(A+BK)$.
Indeed,
\bea \label{almostRiccati}
-S + A'S(A+BK)+ K'B'S(A+BK) = -Q+K'B'P(A+BK).
\eea
We know that 
\bea \label{itsRiccati}
-S + A'S(A+BK)= -Q,
\eea
is the discrete-time algebraic Riccati equation (DTARE).
Subtracting (\ref{itsRiccati}) from (\ref{almostRiccati}), we have
\beann
K'B'S(A+BK)&=&K'B'P(A+BK).
\eeann
Thus,
\bea \label{sisp}
S=P
\eea
and $S$ satisfies the DTARE.  Therefore, we have a diagonalized system

\bea \label{nondiffham3} 
\left[ \begin{array}{c}
x^+ \\
z \\
\end{array}\right]&=&
\left[\begin{array}{cc}
A+BK & 0\\
0    & (A+BK)'\\
\end{array}\right]
\left[\begin{array}{c}
x \\
z^+ \\
\end{array}\right]+
\left[\begin{array}{c}
f_{\psi}(x,x^+)\\
h_{\psi}(x,x^+)\\
\end{array}\right]
\eea
where 
\beann
f_{\psi}(x,x^+)&=&(I+BR^{-1}B'P)^{-1}(F(x,\psi(x^+)) - BR^{-1}B'\psi(x^+))
\eeann
and
\beann
h_{\psi}(x,x^+)=(A+BK)'Pf_{\psi}(x,x^+) + g_{\psi}(x,x^+).~~~~~~~~~~~~~~~
\eeann
The nonlinear terms $f_{\psi}$ and $h_{\psi}$ are $\mathcal{C}^r$ functions for $r \geq 1$ such that
\bea 
f_{\psi}(0,0)=0, && h_{\psi}(0,0)=0 \label{fhcond1}\\
\frac{\partial f_{\psi}}{\partial(x,x^+)}(0,0)=0, && \frac{\partial h_{\psi}}{\partial (x,x^+)}(0,0)=0.\label{fhcond2}
\eea
because of (\ref{fgcond1}) and (\ref{fgcond2}).

\subsection{The Local Stable Manifold Theorem}

Given the original dynamics (\ref{nondiffham}) we look for the local stable manifold described by $\lambda=\phi(x)$.  Since the linear term of the stable manifold for the system (\ref{nondiffham}) is $Px$ where $P$ is the solution to DTARE (\ref{dtare}), we assume the local stable manifold is of the form
\bea \label{transform1}
\lambda=Px+ \psi(x)
\eea
where $\psi(x)$ only contains the nonlinear term of $\phi(x)$.  In the two-step process of diagonalization of the system (\ref{nondiffham}), we introduce the $z$ coordinate through the transformation
\[\lambda=z+Sx.\]
Since $S=P$, it must be that
\beann
z= \lambda - Px = Px + \psi(x) -Px = \psi(x).
\eeann
Then, it suffices to prove existence of the local stable manifold $z=\psi(x)$ for the diagonalized system (\ref{nondiffham3}).  In order to show the existence of $z=\psi(x)$, we use the Contraction Mapping Principle (CMP).  To invoke the CMP, we will need a map $T:\mathbb{X}\longrightarrow \mathbb{X}$ that is a contraction on a complete metric space $\mathbb{X}$.

\begin{theorem} \label{nondiffthm}
Given the dynamics in (\ref{nondiffham3}) with the nonlinear terms $f_{\psi}$ and $h_{\psi}$ that are $\mathcal{C}^r$ functions satisfying the conditions (\ref{fhcond1}) and (\ref{fhcond2}) and a hyperbolic fixed point $0 \in \mathbb{R}^{2n}$, there exists a local stable manifold
\bea \label{psimanifold}
z=\psi(x)
\eea
around the fixed point $0$ where $\psi$ is a $\mathcal{C}^r$ function. 
\end{theorem}

\begin{proof}

First notice that $f_{\psi}$ and $g_{\psi}$ are cut-off functions.  It follows that $h_{\psi}$ is also a cut-off function.  It suffices to prove the theorem for some $\epsilon > 0$ since the cut-off functions $f_{\psi}(x,x^+;\epsilon)$ and $h_{\psi}(x,x^+;\epsilon)$ agree with $f_{\psi}(x,x^+)$ and $h_{\psi}(x,x^+)$ for $\vt x \vt,~\vt x^+ \vt \leq \epsilon$.

By (\ref{fbnd1})-(\ref{gbnd2}), (\ref{cond1})-(\ref{cond2}), and (\ref{fhcond1})-(\ref{fhcond2}), there exists $\bar{N}_1,~\bar{N}_2 >0$ such that 
\bea
\vt f_{\psi}(x,y;\epsilon)-f_{\psi}(\tilde{x},\tilde{y};\epsilon) \vt &\leq& \bar{N}_1 \epsilon \Big[ \vert x- \tilde{x} \vert + \vt  y -\tilde{y} \vt \Big] \label{fphibnd1}\\
\vt h_{\psi}(x,y;\epsilon)-h_{\psi}(\tilde{x},\tilde{y};\epsilon) \vt &\leq& \bar{N}_1 \epsilon \Big[ \vert x- \tilde{x} \vert + \vt  y -\tilde{y} \vt \Big] \label{gphibnd1}
\eea
and 
\bea
\Big \vt \frac{\partial f_{\psi}}{\partial (x,y)}(x,y;\epsilon)- \frac{\partial f_{\psi}}{\partial (x,y)} (\tilde{x},\tilde{y};\epsilon) \Big \vt &\leq&  \bar{N}_2 \Big[ \vert x- \tilde{x} \vert + \vt  y -\tilde{y} \vt \Big] \label{fphibnd2}\\
\Big \vt \frac{\partial h_{\psi}}{\partial (x,y)}(x,y;\epsilon)-\frac{\partial h_{\psi}}{\partial (x,y)}(\tilde{x},\tilde{y};\epsilon) \Big \vt &\leq& \bar{N}_2 \Big[ \vert x- \tilde{x} \vert + \vt  y -\tilde{y} \vt \Big]. \label{gphibnd2}
\eea

Moreover, from the bound (\ref{origmvt1}) we know the following
\bea \label{mvtfphi}
\Big \vt  \frac{\partial f_{\psi}}{\partial (x,y)}(\xi_1,\xi_2;\epsilon) \Big \vt \leq \bar{N}\epsilon
\eea
and
\bea \label{mvtgphi}
\Big \vt  \frac{\partial h_{\psi}}{\partial (x,y)}(\xi_1,\xi_2;\epsilon) \Big \vt \leq \bar{N} \epsilon.
\eea
for $\bar{N}>0$ and $\vt \xi_1 \vt, \vt \xi_2 \vt < \epsilon$.

Henceforth we suppress $\epsilon$ and write $f_{\psi}(x,y), h_{\psi}(x,y)$ for $f_{\psi}(x,y;\epsilon),~h_{\psi}(x,y;\epsilon)$.

Let $l(\epsilon)$ with $l(0)=0$ and $\psi \in \mathbb{X}\subset \mathcal{C}^0({\vt x \vt \leq \epsilon})$ where $\mathbb{X}$ is space of $\psi: \mathbb{R}^n \longrightarrow \mathbb{R}^n$ such that
\bea
\psi(0)&=&0 \label{Xcond1}\\
\vt \psi(x)-\psi(\bar{x}) \vt &\leq& l(\epsilon)\vt x-\bar{x}  \vt \label{Xcond2}
\eea
for $x,~\bar{x} \in \bar{\mathcal{B}}_{\epsilon}(0) \subset \mathbb{R}^n$ where $\bar{\mathcal{B}}_{\epsilon}(0)$ is closed ball around $0$ with radius $\epsilon$.
We define 
\bea \label{thenorm}
\V \psi \V= \sup_{\vt x \vt \leq \epsilon}\left\vt \frac{\psi(x)}{x} \right\vt.
\eea

To show $\mathbb{X}$ is a complete metric space it suffices to show that $\mathbb{X}$ is closed since $\mathbb{X} \subset \mathcal{C}^0(\{\vt x \vt \leq \epsilon\})$.  We take a sequence $\{\psi_n \} \in 
\mathbb{X}$ such that $\psi_n  \longrightarrow  \psi$ in $\mathcal{C}^0$ norm. For large $N>n$, $\vert \psi_n(x)-\psi(x) \vert 
\leq \frac{\epsilon}{2}$ for all $x \in \bar{\mathcal{B}}_{\epsilon}(0)$.  Then,
\beann
\left\vert \psi(x) - \psi(\bar{x})  \right\vert &\leq& \left\vt \psi(x) -\psi_n(x)  \right\vt + \left\vt \psi_n(x) - \psi_n(\bar{x}) \right\vt + \left\vt \psi_n(\bar{x}) - \psi(\bar{x}) \right\vt\\
& \leq & \frac{\epsilon}{2} + l(\epsilon) \vt x-\bar{x} \vt  +\frac{\epsilon}{2}.
\eeann
By letting $\epsilon \rightarrow 0$, we have that $\vert \psi(x) - \psi(\bar{x}) \vert \leq l(\epsilon) \vt x-\bar{x} \vt$.  Thus, $\psi $ is a Lipschitz function.  Similarly, the condition (\ref{Xcond1}) is easily satisfied.  It follows that
\beann
\vt \psi(0)-0 \vt &\leq& \vt \psi(0)- \psi_n(0) \vt + \vt \psi_n(0) - 0 \vt \leq \epsilon.
\eeann
Thus, $\psi \in \mathbb{X}$.  Hence $\mathbb{X}$ is closed.  
Moreover, $\mathbb{X}$ is a complete metric space with the norm defined on (\ref{thenorm}).

Solving the $z$ dynamics in (\ref{nondiffham3}) via the variation of constants formula, we have
\bea \label{nlinearlamsol}
z_j&=& (A'+K'B')^{k-j}z_k + \sum_{l=j}^{k-1}(A'+K'B')^{l-j}h_{\psi}(x_l,x_{l+1})
\eea
for $j<k$.  Let $j=0$ and $k=\infty$, then (\ref{nlinearlamsol}) changes to
\bea \label{nlinearlamsol2}
z_0&=& \sum_{l=0}^{\infty}(A'+K'B')^{l}h_{\psi}(x_l,x_{l+1})
\eea

We define a mapping $T: \mathbb{X} \longrightarrow \mathbb{X}$ by 
\bea \label{Tmap}
(T\psi)(x_0)= \sum_{l=0}^{\infty}(A'+K'B')^{l}h_{\psi}(x_l,x_{l+1}).
\eea
where $x_l, x_{l+1}$ satisfy
\[x^+ = A+BK x + f_{\psi}(x,x^+).\]

From this fixed point equation, we look for the solution
\beann
(T\psi)(x_0)=\psi(x_0).
\eeann
We must show $T\psi \in \mathbb{X}$ and prove $T$ is a contraction on $\mathbb{X}$.


Suppose $x_0=0$ is the initial condition. Clearly from the equation (\ref{nlinearxmap2})-(\ref{consxeq3}), $x_k=0$ for all $k$. Together with $x_k=0$ for all $k$ and the condition (\ref{fhcond2}) $h_{\psi}(0,0)=0$, we have 
\bea \label{sXcond1}
T\psi(0)=0.
\eea
Hence $T\psi$ satisfies the condition (\ref{Xcond1}).

We now prove the Lipschitz condition (\ref{Xcond2}) for $T\psi$.

For $\psi \in \mathbb{X}$ and the initial conditions $x_0,~\bar{x}_0 \in \mathbb{R}^n$, we denote $x_k=x(k,x_0,\psi)$
to be the solution of the state dynamics.  
\beann 
x^+&=&(A+BK)x + f_{\psi}(x,x^+)\\
x(0)&=&x_0. 
\eeann
Similarly, for $\psi \in \mathbb{X}$ and the initial conditions $\bar{x}_0 \in \mathbb{R}^n$, let $x_k=x(k,\bar{x}_0,\psi)$ the solution of
\beann 
x^+&=&(A+BK)x + f_{\psi}(x,x^+)\\
x(0)&=&\bar{x}_0. 
\eeann
Recall $\alpha=\max_j \vt \lambda_j \vt$ where $\lambda_j \in \sigma(A+BK)$ and $\vt \lambda_j \vt < 1$. Using the estimate (\ref{fphibnd1}), at one-time step 
\beann
\vt x_{k+1}- \bar{x}_{k+1}\vt \leq \alpha \vt x_k - \bar{x}_k \vt + \bar{N}_1 \left[ \vt x_k-\bar{x}_k  \vt + \vt x_{k+1}- \bar{x}_{k+1}\vt   \right].
\eeann
For $1-\bar{N}_1\epsilon > 0$, 
\beann
\vt x_{k+1}- \bar{x}_{k+1}\vt \leq \frac{\alpha + \bar{N}_1\epsilon}{1-\bar{N}_1\epsilon} \vt x_k - \bar{x}_k \vt
\eeann
and recursively,
\beann
\vt x_{k}- \bar{x}_{k}\vt \leq \left( \frac{\alpha + \bar{N}_1\epsilon}{1-\bar{N}_1\epsilon} \right)^k \vt x_0 - \bar{x}_0 \vt
\eeann
As long as
\[\epsilon < \frac{1-\alpha}{2\bar{N}_1},\]
then
\[ \frac{\alpha + \bar{N}_1\epsilon}{1-\bar{N}_1\epsilon} < 1.\]
Thus,
\bea \label{uniformx1}
\vt x_k - \bar{x}_k  \vt \leq \vt x_0 - \bar{x}_0 \vt.  
\eea

Using the bounds (\ref{gphibnd1}) and (\ref{uniformx1}),
\beann
\vt T\psi(x_0)- T\psi(\bar{x}_0)\vt &\leq& \left\vt  \sum_{l=0}^{\infty}(A'+K'B')^{l}( h_{\psi}(x_l,x_{l+1})-h_{\psi}(\bar{x}_l,\bar{x}_{l+1}) )\right\vt\\
&\leq&   \sum_{l=0}^{\infty}\alpha^{l}\left\vt  h_{\psi}(x_l,x_{l+1})-h_{\psi}(\bar{x}_l,\bar{x}_{l+1}) \right\vt \\ 
&\leq&   \sum_{l=0}^{\infty}\alpha^{l} \bar{N}_1\epsilon \left[ \vt x_l - \bar{x}_l \vt + \vt x_{l+1} - \bar{x}_{l+1} \vt \right]\\
&\leq&   \sum_{l=0}^{\infty}\alpha^{l} 2\bar{N}_1\epsilon \vt x_0-\bar{x}_0 \vt \\
&\leq& \frac{2\bar{N}_1\epsilon}{1-\alpha} \vt x_0 - \bar{x}_0 \vt.
\eeann
Let
\[l(\epsilon)= \frac{2\bar{N}_1\epsilon}{1-\alpha}.\]
Notice that $l(\epsilon) \rightarrow 0$ as $\epsilon \rightarrow 0$.
Thus,
\beann
\vt T\psi(x_0)- T\psi(\bar{x}_0)\vt &\leq& l(\epsilon) \vt  x_0 - \bar{x}_0  \vt
\eeann
and so $(T\psi)$ satisfies the condition (\ref{Xcond2}) for $\epsilon >0$ sufficiently small.

Hence $T$ maps from $\mathbb{X} \rightarrow \mathbb{X}.$


Next we show T is a contraction on $\mathbb{X}$.

We express the solutions to the state dynamics
\[x_k=x(k,x_0,\psi)\]
and
\[\bar{x}_k=\bar{x}(k,\bar{x}_0,\psi)\]
in the implicit form,
\beann
x_k=(A+BK)^k x_0 + \sum_{j=0}^{k-1} (A+BK)^{k-1-j}f_{\psi}(x_j,x_{j+1})
\eeann
and
\beann
\bar{x}_k=(A+BK)^k \bar{x}_0 + \sum_{j=0}^{k-1} (A+BK)^{k-1-j}f_{\psi}(\bar{x}_j,\bar{x}_{j+1}),
\eeann
respectively.

We now denote $x_j=x(j,x_0,\psi)$ and $\bar{x}_j=\bar{x}(j,x_0,\bar{\psi})$ be the solutions to the state dynamics and satisfy the implicit form equations
\[x_k=(A+BK)^k x_0 + \sum_{j=0}^{k-1} (A+BK)^{k-1-j}f_{\psi}(x_j,x_{j+1})\]
and
\[\bar{x}_k=(A+BK)^k x_0 + \sum_{j=0}^{k-1} (A+BK)^{k-1-j}f_{\bar{\psi}}(\bar{x}_j,\bar{x}_{j+1}),\]
respectively.  

The estimates (\ref{fphibnd1})-(\ref{gphibnd1}) with the trajectories $x(j,x_0,\psi)$ and $x(j,x_0,\bar{\psi})$ becomes
\bea
\vt f_{\psi}(x,y)-f_{\psi}(\tilde{x},\tilde{y}) \vt &\leq& r_1(\epsilon)\vt y-\bar{y} \vt + r_2(\epsilon) \V \psi- \bar{\psi} \V + r_3(\epsilon)\vt x-\bar{x} \vt \label{fphibnd2}\\
\vt h_{\psi}(x,y)-h_{\psi}(\tilde{x},\tilde{y}) \vt &\leq& r_1(\epsilon)\vt y-\bar{y} \vt + r_2(\epsilon) \V \psi- \bar{\psi} \V + r_3(\epsilon)\vt x-\bar{x} \vt  \label{gphibnd2}
\eea 
where
\beann
r_1(\epsilon)&=& n_{1,1}l(\epsilon) + n_{1,2}\epsilon + n_{1,3}l(\epsilon)\epsilon\\
r_2(\epsilon)&=& n_{2,1}\epsilon + n_{2,2}\epsilon^2\\
r_3(\epsilon)&=& n_3\epsilon
\eeann
and $n_{i,j}$ are positive constants.  Observe that $r_i(\epsilon) \rightarrow 0$ as $\epsilon \rightarrow 0$.

At one-time step,
\beann
\vt x_{k+1}-\bar{x}_{k+1} \vt &\leq& m_2(\epsilon) \vt x_k-\bar{x}_k \vt + m_3(\epsilon) \V \psi-\bar{\psi} \V
\eeann
where
\[m_2(\epsilon)=\frac{\alpha + r_3(\epsilon)}{1-r_1(\epsilon)} \]
and
\[m_3(\epsilon)=\frac{r_2(\epsilon)}{1-r_1(\epsilon)}.\]
By invoking Gronwall's inequality (\ref{gronwall1}) and assuming that for some small $\epsilon > 0$
\[m_2(\epsilon)<1,\]
then
\bea \label{xphibound}
\vt x_{k}-\bar{x}_{k} \vt &\leq& m_3(\epsilon)\left[ \sum_{j=0}^{k-1} m_2(\epsilon)^{k-1-j}\right] \V \psi-\bar{\psi} \V \nn\\
&\leq& m_3(\epsilon) \frac{1-m_2(\epsilon)^k}{1-m_2(\epsilon)} \V \psi-\bar{\psi} \V \nn \\
&\leq& \frac{m_3(\epsilon)}{1-m_2(\epsilon)} \V \psi-\bar{\psi} \V 
\eea

With the bounds (\ref{gphibnd2}) and (\ref{xphibound}), we get
\beann
\vt T\psi(x_0)- T\bar{\psi}(x_0)\vt &\leq& \left\vt  \sum_{l=0}^{\infty}(A'+K'B')^{l}( h_{\psi}(x_l,x_{l+1})-h_{\bar{\psi}}(\bar{x}_l,\bar{x}_{l+1}) )\right\vt\\
&\leq&  \sum_{l=0}^{\infty}\alpha^{l} \left\vt h_{\psi}(x_l,x_{l+1})-h_{\bar{\psi}}(\bar{x}_l,\bar{x}_{l+1}) \right\vt \\
&\leq&  \sum_{l=0}^{\infty}\alpha^{l} \left[ r_1(\epsilon)\vt x_{l+1}-\bar{x}_{l+1} \vt + r_2(\epsilon) \V \psi- \bar{\psi} \V + r_3(\epsilon)\vt x_l-\bar{x}_l \vt \right]\\
&\leq&  \sum_{l=0}^{\infty}\alpha^{l} \left[ 2(r_1(\epsilon)+ r_3(\epsilon)) \frac{m_3(\epsilon)}{1-m_2(\epsilon)} \V \psi-\bar{\psi} \V + r_2(\epsilon) \V \psi- \bar{\psi} \V  \right].
\eeann

It follows that
\beann
\left\vt T\psi(x_0)- T\bar{\psi}(\bar{x}_0) \right\vt &\leq& c(\epsilon) \V \psi - \bar{\psi} \V
\eeann
where
\[c(\epsilon)=\frac{2m_3(\epsilon)(r_1(\epsilon)+ r_3(\epsilon))}{(1-m_2(\epsilon))(1-\alpha)} + \frac{r_2(\epsilon) }{1-\alpha}.\]
Notice that $c(\epsilon) \rightarrow 0$ as $\epsilon \rightarrow 0$.
Thus, for sufficiently small $\epsilon>0$
\[c(\epsilon)<1.\]
Hence, $T$ is a contraction on $\mathbb{X}$ for $\epsilon$ sufficiently small.  Therefore, there exists a unique $\psi \in \mathbb{X}$ such that
\[\psi=T\psi.\]

Let $\bar{x},\bar{z}$ satisfy
\[x^+= A+BK x + f_{\phi}(x,x^+)\]
and
\[z^+= (A+BK)'z + g_{\phi}(x,x^+),\]
respectively.  Choose $\bar{z}=\psi(\bar{x})$.  Using the definition of the contraction T in (\ref{Tmap}), we have
\beann
\bar{z}_0&=&(T\psi)(\bar{x}_0)\\
&=& \sum_{l=0}^{\infty}(A'+K'B')^{l}h_{\psi}(\bar{x}_l,{x}_{l+1}).
\eeann
It follows that $(\bar{x},\bar{z})$ is a solution to the difference equation (\ref{nondiffham3}) and (\ref{psimanifold}) defines a $\mathcal{C}^1$ invariant manifold.  Thus, (\ref{phimanifold}) is the graph of a $\mathcal{C}^1$ invariant manifold.
As described in \cite{Kr01}, one can show straight forwardly that (\ref{phimanifold}) defines $\mathcal{C}^r$ invariant manifold given the data are $\mathcal{C}^r$ smooth.
\end{proof}

\section{Some Properties}

\subsection{Eigenstructure}

Recall from Section $1$ the bidirectional linear Hamiltonian dynamics
\bea \label{eigdyn}
\left[ \begin{array}{c}
x^+ \\
\lambda \\
\end{array}\right]&=&
\mathbb{H}
\left[\begin{array}{c}
x \\
\lambda^+ \\
\end{array}\right]
\eea
where
\beann
\mathbb{H}=\left[\begin{array}{cc}
A & -BR^{-1}B'\\
Q & A'\\
\end{array}\right].
\eeann

\begin{definition} \label{sfactor}
Suppose
\bea \label{therelation}
\mathbb{H}
\left(\begin{array}{c}
\delta x \\
\mu \delta \lambda \\
\end{array}\right)=
\left(\begin{array}{c}
 \mu \delta x \\
\delta \lambda \\
\end{array}\right).
\eea
Then we call $\mu$ the eigenvalue of the dynamics (\ref{eigdyn}).
\end{definition}

We would like to show that the linear bidirectional Hamiltonian matrix $\mathbb{H}$ is hyperbolic; i.e., the eigenvalue of the dynamics (\ref{eigdyn}) lies strictly inside and outside the unit circle.

\begin{theorem}
If $\mu$ is an eigenvalue of the dynamics (\ref{eigdyn}) satisfying the relation (\ref{therelation}), then $\frac{1}{\mu}$ is also an eigenvalue of the dynamics (\ref{eigdyn}).
\end{theorem}
\begin{proof}
First, we decompose $\mathbb{H}$ into
\beann
\left[\begin{array}{cc}
0 & I\\
I & 0\\
\end{array}\right]
\left[\begin{array}{cc}
Q & A'\\
A & -BR^{-1}B'\\
\end{array}\right]
=\mathbb{H}.
\eeann
Let's call
\beann
\mathbb{S}=\left[\begin{array}{cc}
Q & A'\\
A & -BR^{-1}B'\\
\end{array}\right],~
\mathbb{I}_{\mu}=\left[\begin{array}{cc}
I & 0\\
0 & \mu I\\
\end{array}\right],~
\mathbb{I}^{\mu}=\left[\begin{array}{cc}
\mu I & 0\\
0    & I\\
\end{array}\right]
\eeann
and notice that $\mathbb{S}$ is symmetric.

We rewrite (\ref{therelation}) to
\beann
\mathbb{H}~
\mathbb{I}_{\mu}
\left(\begin{array}{c}
\delta x \\
\delta \lambda \\
\end{array}\right)
=
\mathbb{I}^{\mu}
\left(\begin{array}{c}
\delta x \\
\delta \lambda \\
\end{array}\right).
\eeann
It follows that
\bea \label{rel2}
\mathbb{I}^{\frac{1}{\mu}}~
\mathbb{H}~
\mathbb{I}_{\mu}
\left(\begin{array}{c}
\delta x \\
\delta \lambda \\
\end{array}\right)
=
\left(\begin{array}{c}
\delta x \\
\delta \lambda \\
\end{array}\right).
\eea
We denote
\beann
\mathbb{H}_{\mu}=
\mathbb{I}^{\frac{1}{\mu}}~
\mathbb{H}~
\mathbb{I}_{\mu}.
\eeann
Observe that $1$ is an eigenvalue of $\mathbb{H}_{\mu}$ from (\ref{rel2}).  Then,
\[det \left[I- \mathbb{H}_{\mu} \right]=0 \Longrightarrow det \left[I- \mathbb{H}'_{\mu} \right]=0\]
where
\beann
\mathbb{H}'_{\mu}=
\mathbb{I}_{\mu}
\mathbb{S}
\left(\begin{array}{cc}
0 & I \\
I & 0 \\
\end{array}\right)
\mathbb{I}^{\frac{1}{\mu}}.
\eeann
Again, $1$ is an eigenvalue of $\mathbb{H}_{\mu}'$; i.e.,
\beann \label{theadjoint}
\mathbb{H}'_{\mu}
\left(\begin{array}{c}
\widetilde{\delta x}\\
\widetilde{\delta \lambda} \\
\end{array}\right)=
\left(\begin{array}{c}
\widetilde{\delta x}\\
\widetilde{\delta \lambda} \\
\end{array}\right).
\eeann
By multiplying $\mathbb{I}_{\frac{1}{\mu}}$ to equation above, we have
\beann
\mathbb{S}
\left(\begin{array}{cc}
0 & I \\
I & 0 \\
\end{array}\right)
\left(\begin{array}{c}
\frac{1}{\mu} \widetilde{\delta x}\\
\widetilde{\delta \lambda}\\
\end{array}\right)=
\left(\begin{array}{c}
\widetilde{\delta x}\\
\frac{1}{\mu} \widetilde{\delta \lambda}\\
\end{array}\right)
\Longleftrightarrow
\mathbb{S}
\left(\begin{array}{c}
\widetilde{\delta \lambda}\\
\frac{1}{\mu} \widetilde{\delta x}\\
\end{array}\right)=
\left(\begin{array}{cc}
0 & I \\
I & 0 \\
\end{array}\right)
\left(\begin{array}{c}
\frac{1}{\mu} \widetilde{\delta \lambda}\\
\widetilde{\delta x}\\
\end{array}\right).
\eeann
Thus,
\beann \label{theadjoint}
\mathbb{H}
\left(\begin{array}{c}
\widetilde{\delta \lambda} \\
\frac{1}{\mu} \widetilde{\delta x}\\
\end{array}\right)=
\left(\begin{array}{c}
\frac{1}{\mu} \widetilde{\delta \lambda} \\
\widetilde{\delta x}\\
\end{array}\right).
\eeann
Hence, $\frac{1}{\mu}$ is an eigenvalue of (\ref{eigdyn}).
\end{proof}
Note that the dynamics (\ref{eigdyn}) admits the infinite eigenvalues, $0$ and $\infty$, due to the singularity of $A$.  

\subsection{Symplectic Form}
The nonlinear dynamics tangent to (\ref{nondiffham}) as derived explicitly using perturbation technique in cite{Na02} is
\bea \label{nonhampert2}
\left[ \begin{array}{c}
\delta x^+\\
\delta \lambda \\
\end{array}\right]=
\left[ \begin{array}{cc}
H_{\lambda ^+ x} & H_{\lambda^+ \lambda^+}\\
H_{x x}          & H_{\lambda^+ x}\\
\end{array}\right]
\left[ \begin{array}{c}
\delta x\\
\delta \lambda^+\\
\end{array}\right]
\eea
where $H_{\lambda ^+ x}, H_{x x}$, and $H_{\lambda^+ \lambda^+}$ are the partial derivatives of the Hamiltonian defined in (\ref{theham}) and $(\delta x, \delta \lambda)$ are tangent vectors in $T_{(x,\lambda)}\mathcal{M}$ for $\mathcal{M}=\{(x,\lambda)| x \in \mathbb{R}^n,~\lambda \in \mathbb{R}^n \}$.
The nondegenerate and bilinear symplectic two-form $\Omega:T_{(x,\lambda)}\mathcal{M} \times T_{(x,\lambda)}\mathcal{M} \mapsto \mathbb{R}$ is 
\beann
\Omega(v,w)=v'Jw
\eeann
with
\[\Omega(v,w)=-\Omega(w,v)\]
where the symplectic matrix,
\beann
J=
\left[ \begin{array}{cc}
0  & I \\
-I & 0 \\
\end{array}\right],
\eeann
and
\[ v=\left[ \begin{array}{c}
\delta x\\
\delta \lambda \\
\end{array}\right],~~~ w=\left[ \begin{array}{c}
\widetilde{\delta x}\\
\widetilde{\delta \lambda} \\
\end{array}\right] \in T_{(x,\lambda)}\mathcal{M}.\]

We would like to show that under the tangent dynamics (\ref{nonhampert2}), the two-form $\Omega$ is invariant; i.e.,
\bea \label{preservingform}
\Omega(v,w)=\Omega(v^+,w^+).
\eea
Then,
\beann \label{Omegalhs}
\Omega(v,w)&=&v'Jw \nn \\
&=&-\delta \lambda '\widetilde{\delta x} + \delta x' \widetilde{\delta \lambda}.
\eeann
By substituting the tangent dynamics of $\delta x$ in (\ref{nonhampert2}) into the above equation, we have
\bea \label{Omegalhs}
\Omega(v,w)&=&-(H_{xx} \delta x + H_{\lambda ^+ x} \delta \lambda ^+ )' \widetilde{\delta x} + \delta x'(H_{xx}\widetilde{\delta x} + H_{\lambda ^+ x}\widetilde{\delta \lambda ^+} )\nn \\
&=&-\delta \lambda ^+ H'_{\lambda ^+ x} \widetilde{\delta x} + \delta x H_{\lambda ^+ x} \widetilde{\delta \lambda ^+}. 
\eea

Similarly,
\bea \label{Omegarhs}
\Omega(v^+,w^+)&=&v^{+'}Jw^+ \nn \\
&=& - (\delta \lambda ^+)' \widetilde{\delta x ^+} + (\delta x^+)' \widetilde{\delta \lambda ^+} \nn \\
&=& - (\delta \lambda ^+)' (H_{\lambda ^+ x} \widetilde{\delta x} + H_{\lambda ^+ \lambda ^+} \widetilde{\delta \lambda ^+} ) + (H_{\lambda ^+ x} \delta x + H_{\lambda ^+ \lambda ^+} \delta \lambda ^+ )' \widetilde{\delta \lambda ^+} \nn \\
&=& -\delta \lambda ^+ H'_{\lambda ^+ x} \widetilde{\delta x} + \delta x H_{\lambda ^+ x} \widetilde{\delta \lambda ^+}.
\eea
Since (\ref{Omegarhs}) and (\ref{Omegalhs}) are equal, then
\[\Omega(v,w)=\Omega(v^+,w^+).\]
Thus, for any two tangent vectors satisfying the dynamics (\ref{nonhampert2}), the value of $\Omega$ does not change.

\section{Lagrangian Submanifold}

The two-form calculated from the last section is 
\beann \label{Omega2form}
\Omega(v,w)=-\delta \lambda^{+'} H_{\lambda^+ x}\widetilde{\delta x} + \delta x' H_{\lambda^+ x} \widetilde{\delta \lambda}^+.
\eeann 
We calculate the  state dynamics tangent to (\ref{nondiffham3}) around the trajectories $(x_j,x_{j+1})$ as
\bea \label{deltaxtandynamics}
\delta x_{j+1}   &=& \Big( (A+BK)  + \frac{\partial f_{\psi}}{\partial x_j}(x_j,x_{j+1}) \Big) \delta x_j \label{tandyn1}\\
&&~~~~~~~~~~~+ \frac{\partial f_{\psi}}{\partial x_{j+1}}(x_j, x_{j+1}) \delta x_{j+1} \nn 
\eea

By the Inverse Function Theorem, we can choose $\epsilon >0 $ small enough so that $I-\frac{\partial f_{\psi}}{\partial x_{k+1}}(x_k, x_{k+1})$ is invertible for $ \vt x_k \vt, \vt x_{k+1} \vt < \epsilon$ because of (\ref{fgcond1}) and hence,
\[I-\frac{\partial f_{\psi}}{\partial x_{k+1}}(0, 0)= I.\]

It follows that (\ref{deltaxtandynamics}) is equivalent to
\beann
\delta x_{k+1}= \Big( I-\frac{\partial f_{\psi}}{\partial x_{k+1}}(x_k, x_{k+1}) \Big)^{-1} \Big( (A+BK)  + \frac{\partial f_{\psi}}{\partial x_k}(x_k,x_{k+1})\Big) \delta x_k
\eeann
and
\bea \label{deltaxsoln}
\delta x_{k+1}&=& \prod_{i=0}^{k}\Big( I-\frac{\partial f_{\psi}}{\partial x_{i+1}}(x_i, x_{i+1})\Big)^{-1} \cdot \\
&&~~~~~~~~~~\Big( (A+BK)  + \frac{\partial f_{\psi}}{\partial x_i}(x_i,x_{i+1})\Big) \delta x_0. \nn
\eea
for $ \vt x_k \vt, \vt x_{k+1} \vt < \epsilon$.  As we let $k \rightarrow \infty$, we have $x_k \rightarrow 0$, $\frac{\partial f_{\psi}}{\partial (x,x^+)}(0) \rightarrow 0$.  It follows that from (\ref{deltaxsoln})
\[\delta x_{k+1} \rightarrow (A+BK)^k \delta x_0 \rightarrow 0.\]

Thus, 
\[\Omega(v,w)= -\delta \lambda^{+'} H_{\lambda^+ x}\widetilde{\delta x} + \delta x' H_{\lambda^+ x} \widetilde{\delta \lambda}^+ \rightarrow 0\]
since $\delta x_k \rightarrow 0$ as $k \rightarrow \infty$ for $v,w$ restricted to the tangent dynamics.  Hence, the local stable manifold is a Lagrangian submanifold.

Denote $W_s$ as the local stable manifold described by the graph $\lambda = \phi(x)$.  The basis for $T_{(x,\phi(x))}W_s$ is of the form
\beann
\frac{\partial}{\partial x_i}
\left[ \begin{array}{c}
x_1\\
\vdots\\
x_{n}\\
\phi_1\\
\vdots\\
\phi_n\\
\end{array}\right]\;\;\mbox{where}\;i=1,\ldots,n.
\eeann
Then, the two-form
\[\Omega(v^i,w^j)=(v^i)' J w^j\]
where
\beann
v^i=
\frac{\partial}{\partial x_{i}}
\left[ \begin{array}{c}
x_{1}\\
\vdots\\
x_{n}\\
\phi_1\\
\vdots\\
\phi_n\\
\end{array}\right],\;w^j=
\frac{\partial}{\partial x_{j}}
\left[ \begin{array}{c}
x_{1}\\
\vdots\\
x_{n}\\
\phi_1\\
\vdots\\
\phi_n\\
\end{array}\right] \in T_{(x,\phi(x))}W_s
\eeann
is equivalent to
\beann \label{closedtoo}
\frac{\partial \phi_i}{\partial x_{j}}-\frac{\partial \phi_j}{\partial
x_{i}}=0\;\;\mbox{  for  }i,j=1,\ldots,n.
\eeann

The equation above implies that $\phi(x)$ is closed.  Then, by 
\emph{Stokes' Theorem} there exists $\pi \in \mathcal{C}^{r}(\mathbb{R}^n)$
such that
\bea \label{voila}
\phi(x)=\frac{\partial \pi}{\partial x}(x)\mbox{  where  }\phi(0)=0
\eea
locally on some neighborhood of 0.  Thus, there exists $\pi \in \mathcal{C}^{r}$ such that
\bea \label{thegradient}
\lambda=\frac{\partial \pi}{\partial x}(x).
\eea

Hence, the local stable manifold $\lambda$ in (\ref{thegradient}) is the gradient of the optimal cost for the bidirectional Hamiltonian dynamics (\ref{theham}).

\section{Application to Dynamic Programming Equations}
Recall from Section $2$ the formulation of a discrete in time infinite horizon optimal control problem of minimizing the cost functional, 
\[ \min_u \sum_{k=0}^{\infty} l(x_k,u_k)\]
subject to the dynamics
\beann
x^+&=&f(x,u)\\
x(0)&=&x_0
\eeann

We assume $l(x,u)$ is convex in $x$ and $u$ so that
\beann
\left[ \begin{array}{cc}
Q   & S \\
S^* & R \\
\end{array}\right] \geq 0 
\eeann
and let $R >0.$
In addition, it is assumed that the pair $(A,B)$ is stabilizable and the pair $(A,Q^{1/2})$ is detectable.

Given that $x(0)=x_0$ the optimal value function $\pi(x)$ is defined by
\[\pi(x_0)=\min_u \sum_{k=0}^{\infty} l(x_k,u_k).\]
This value function satisfies a functional equation, called the dynamic programming equation.  The optimal feedback $\kappa(x)$ is constructed from the dynamic programming equation.  We state the optimality principle:
\begin{theorem}
Discrete-Time Optimality Principle:  
\bea \label{optcond}
\pi(x)=\min_{u} \{\pi(f(x,u)) + l(x,u)\}
\eea
\end{theorem}
\begin{proof}
We have that
\beann
\pi(x_0)&=&\min_{u_0} \{\sum_{k=0}^{\infty} l(x_k,u_k)\}\\
        &=&\min_{u_0} \{l(x_0,u_0) +\sum_{k=1}^{\infty} l(x_k,u_k) \}\\
        &=&\min_{u_0} \{l(x_0,u_0) + \pi(x_1)\}
\eeann
Generalizing the optimality principle at the $k^{th}$-step, we have
\bea \label{dhjb1}
\pi(x)=\min_{u} \{\pi(x^+) +  l(x,u)\}.
\eea
\end{proof}
The optimality equation (\ref{optcond}) is the first equation of the dynamic programming equations.  An optimal policy $u^*=\kappa(x)$ must satisfy
\beann
\pi(x) -  \pi(f(x,u^*)) - l(x,u^*)  = 0
\eeann
if we assume convexity of the LHS of (\ref{dhjb1}).  We can find $u^*$ through
\[\frac{\partial (\pi(x) - \pi(f(x,u)) - l(x,u)}{\partial u}=0\]
which by the chain rule becomes 
\[\frac {\partial \pi}{\partial x}(f(x,u)) \frac{\partial f}{\partial u}(x,u) +  \frac{\partial l}{\partial u}(x,u)=0\]
Thus, $\pi(x)$ and $\kappa(x)$ satisfy these equations, the Dynamic Programming Equations (DPE):
\bea
\pi(x) -  \pi(f(x,u)) - l(x,u) ) &=& 0 \label{dpe1}\\
\frac {\partial \pi}{\partial x}(f(x,u)) \frac{\partial f}{\partial u}(x,u) + \frac{\partial l}{\partial u}(x,u) &=& 0 \label{dpe2}
\eea
To show the existence of the local solutions, $\pi(x)$ and $\kappa(x)$, we use the Pontryagin Maximum Principle (Theorem $2.1$). From the condition (\ref{hode3}) and the consequence of Local Stable Manifold Theorem where $\lambda$ is function of $x$; i.e., $\lambda=\frac{\partial \pi}{\partial x}(x)$, we have
\[u^*=\kappa(x)=\mbox{arg}\min_vH(x,\frac{\partial \pi}{\partial x}(x),v).\]
The above equation is equivalent to
\[\frac{\partial H}{\partial u}(x,\frac{\partial \pi}{\partial x}(x),\kappa(x))=0\]
which is essentially (\ref{dpe2}).  Thus, $\frac{\partial \pi}{\partial x}(x)$ and $\kappa(x)$ solve (\ref{dpe2}). Since 
\beann
\lambda=\frac{\partial H}{\partial x}
\eeann
from the PMP and
\[ \lambda=\frac{\partial \pi}{\partial x}, \]
it follows that 
\bea \label{predpe1}
\frac{\partial \pi}{\partial x}= \frac{\partial \pi}{\partial x}(x)\frac{\partial f}{\partial x}(x,u^*)+ \frac{\partial l}{\partial x}(x,u^*)
\eea
Integrating (\ref{predpe1}) w.r.t. $x$, we get
\[\pi(x)-\pi(f(x,\kappa(x)))-l(x,\kappa(x))=0\]
which is (\ref{dpe1}).  Therefore, $\pi$ and $\kappa$ solve the DPE (\ref{dpe1}, \ref{dpe2}).  Furthermore $\pi \in \mathcal{C}^r$ and $\kappa \in \mathcal{C}^{r-1}$ since $l \in \mathcal{C}^r$ and $f \in \mathcal{C}^{r-1}$ in the Hamiltonian.


\section*{Acknowledgements}
The material discussed on this paper was part of my dissertation at the University of California, Davis.  I would like to thank Prof. Arthur J. Krener, my advisor, for his guidance and support.


\end{document}